\documentclass[12pt]{amsart}
\usepackage[top=30truemm,bottom=30truemm,left=25truemm,right=25truemm]{geometry}
\usepackage{txfonts}
\usepackage{mathrsfs}

\usepackage{color}
\usepackage{bm}
\usepackage{amsfonts,amssymb}
\usepackage{dsfont}
\usepackage{amscd}
\usepackage{extarrows}
\usepackage{amsmath}
\usepackage{mathrsfs}
\usepackage{enumerate}
\usepackage{amscd}
\usepackage[all]{xy}
\usepackage[pagebackref,colorlinks]{hyperref}

\newtheorem{thm}{Theorem}[section]
\newtheorem{defn}[thm]{Definition}
\newtheorem{lem}[thm]{Lemma}
\newtheorem{cor}[thm]{Corollary}
\newtheorem{prop}[thm]{Proposition}
\newtheorem{rem}[thm]{Remark}

\newcommand{\be}{\begin{equation}}
	\newcommand{\ee}{\end{equation}}
\newcommand{\bea}{\begin{eqnarray}}
	\newcommand{\eea}{\end{eqnarray}}
\newcommand{\ben}{\begin{eqnarray*}}
	\newcommand{\een}{\end{eqnarray*}}
\newcommand{\bt}{\begin{split}}
	\newcommand{\et}{\end{split}}
\newcommand{\bet}{\begin{equation}}

\begin{document}
	\title{Monge-Amp\`ere type equations on compact Hermitian manifolds}
				
\author[Y. Li]{Yinji Li}
\address{Yinji Li: Institute of Mathematics, Academy of Mathematics and Systems Sciences, Chinese Academy of Sciences,
Beijing, 100190, P. R. China}
\email{1141287853@qq.com}
\author[G. Lin]{Genglong Lin}
\address{Genglong Lin: Institute of Mathematics, Academy of Mathematics and Systems Sciences, Chinese Academy of Sciences,
Beijing, 100190, P. R. China}
\email{lingenglong@amss.ac.cn}
\author[X. Zhou]{Xiangyu Zhou}
\address{ Xiangyu Zhou:\ Institute of Mathematics, Academy of Mathematics and Systems Sciences, and Hua Loo-Keng Key
Laboratory of Mathematics, Chinese Academy of Sciences, Beijing, 100190, P. R. China}
\email{xyzhou@math.ac.cn}
		%\author[F. Deng]{Fusheng Deng}
		%\address{Fusheng Deng: \ School of Mathematical Sciences, University of Chinese Academy of Sciences\\ Beijing 100049, P. R. China}
		%\email{fshdeng@ucas.ac.cn}
		%
		%\author[J. Ning]{Jiafu Ning}
		%\address{Jiafu Ning: \ Department of Mathematics, Central South University, Changsha, Hunan 410083, P. R. China.}
		%\email{jfning@csu.edu.cn}

\begin{abstract}
Given a cohomology $(1,1)$-class $\{\beta\}$ of compact Hermitian manifold $(X,\omega)$ possessing a bounded potential and fixed a model potential $\phi$, motivated by Darvas-Di Nezza-Lu and Li-Wang-Zhou's work, we show that degenerate complex Monge-Amp\`ere equation $(\beta+dd^c \varphi)^n=e^{\lambda \varphi}\mu$ has a unique solution in the relative full mass class $\mathcal{E}(X,\beta,\phi)$, where $\mu$ is a non-pluripolar measure on $X$ and $\lambda\geq0$ is a fixed constant. 

As an application, we give an explicit description of Lelong numbers of elements in $\mathcal{E}(X,\beta,\phi)$ which generalized a theorem of Darvas-Di Nezza-Lu in the Hermitian context.
\end{abstract}

	\maketitle
\noindent{\bf Key words:} Monge-Amp\`ere equation, model potential, Lelong number
\\
\noindent{\bf Mathematics Subject Classification:} 32U05, 32U15, 32U20, 32U25, 32U40
\tableofcontents
\section{Introduction}
Solving degenerate complex Monge-Amp\`ere equations has been the subject of intensive studies in the past decade, in connection with the search for canonical models and metrics of complex algebraic varieties. Many of these results relied on the seminal work of Yau \cite{Yau78} which used a continuity method and a prior estimates to construct smooth solutions to non-degenerate Monge-Amp\`ere equations. It is natrual and necessary to study singular solutions of degenerate equations.

 Later on, based on pluripotential theory on compact Kahler manifolds (see \cite{GZ17}\cite{GZ05}), Guedj-Zeriahi proved that there exists $\varphi\in\mathcal{E}(X,\omega)$ such that 
\begin{equation}\label{eq}
(\omega+dd^c \varphi)^n=\mu
\end{equation}
if and only if $\mu$ does not charge pluripolar sets, where $\omega$ is a K\"ahler form on $X$ and $\mathcal{E}(X,\omega)$ consists of $\omega$-psh functions which has full Monge-Amp\`ere mass on $X$ \cite{GZ07}. Note that the proof of this result depended on Yau's result. In order to look for singular solution on singular spaces, Eyssidieux-Guedj-Zeriahi \cite{EGZ09} used Ko\l odziej's $L^{\infty}$-estimate (see \cite{Kol98}) to obtain bounded solutions  on singular K\"ahler spaces where  RHS of (\ref{eq})  has $L^p(p>1)$ density with respect to smooth volume form. As noticed by \cite{BEGZ10}, non-pluripolar product is well-defined on compact K\"ahler manifolds. Based on a decomposition trick on big classes, they obtained a more general result that there exists $\varphi\in\mathcal{E}(X,\theta)$ such that
\begin{equation} \label{equ1}
(\theta+dd^c \varphi)^n=\mu
\end{equation}
if and only if the probability measure $\mu$ is non-pluripolar, where $\{\theta\}$ is a big cohomology $(1,1)$-class of compact K\"ahler manifold $X$.  

In order to find a proof which does not depend on Yau's result, Berman-Boucksom-Guedj-Zeriahi \cite{BBGZ13} developed variational method and proved the existence and uniqueness of (\ref{equ1}). Their method were generalized by Darvas-Di Nezza-Lu to study Monge-Amp\`ere equations with prescibed singularities \cite{DDNL23}. Berman-Boucksom-Guedj-Zeriahi's method depended on  the regularity of quasi-psh envelope, i.e. $\mbox{MA}(V_{\theta})$ has $L^{\infty}$ density with respect to Lebesgue measure, where $V_{\theta}:=\sup\{u\in\mbox{PSH}(X,\theta),u\leq0\}$. The variational technique also relied on the fact that locally pluripolar sets are $\mbox{PSH}(X,\theta)$-pluripolar.

In \cite{LWZ23} the first, the third author and Zhiwei Wang considered degenerate Monge-Amp\`ere equations with respect to a cohomology class $\{\beta\}$ possessing a bounded potential $\rho$:
$$(\beta+dd^c\varphi)^n=fdV.$$
The right hand side of the equation has $L^p(p>1)$ Lebesgue density. They used Kolodziej-Nguyen's estimate and provided some interesting applications. In \cite{LWZ24} they proved that when the bounded potential $\rho$ is continuous, Capacity $Cap_{\beta}$ can characterize pluripolar sets if and only if Equation (\ref{equ1}) is solvable:
\begin{thm}\cite{LWZ24}
Let $(X,\omega)$ be a compact Hermitian manifold and $\{\beta\}\in H^{(1,1)}(X,\mathbb{R})$ be a real $(1,1)$-class with a smooth representative $\beta$. Suppose that there exists a continuous function $\rho\in PSH(X,\beta)$ then the following two conditions are equivalent:
\begin{itemize}
\item[(1)] For every non-pluripolar Radon measure $\mu$ satisfying $\mu(X)=\int_X \beta^n$, there exists $\varphi\in\mathcal{E}(X,\beta)$ such that $$(\beta+dd^c\varphi)^n=\mu,$$
\item[(2)]$Cap_{\beta}(\cdot)$ can characterize pluripolar sets in the sense that: for every Borel subset $E$, $Cap_{\beta}(E)=0$ if and only if $E$ is pluripolar.
\end{itemize} 
\end{thm}

Our goal is to reformulate Darvas-Di Nezza-Lu's work into a new setting and solve Monge-Amp\`ere equations with prescribed singularity, generalizing the results mentioned as above. Namely, we considered the equation on a cohomology class $\{\beta\}$ possessing a bounded potential $\rho$. Notice that this class is nef and it is conjectured by Tosatti \cite[Conj. 4.1]{Tos23} that nef $(1,1)$-classes on Calabi-Yau manifold always contain closed positive current with bounded potentials. Our key point is to prove the regularity of Monge-Amp\`ere measure of  envelope and equivalence between locally pluripolarity and $\mbox{PSH}(X,\beta)$-polarity under our setting. Then we can use the variational method to solve Monge-Amp\`ere equations if the model potential $\phi$ (see Definition \ref{def of model potential}) has small unbounded locus. Notice that $\rho$ itself is a model potential with small unbounded locus by definition. We will use this special case and follow the line of \cite{DDNL21} where the authors established a supersolution technique. Finally we can drop the restriction on the model potential $\phi$ and prove the following:

\begin{thm}
Suppose that $(X,\omega)$ ia a compact Hermitian manifold and $\{\beta\}$ is a class of $X$ having a bounded potential. Let $\phi$ be a model potential and $\mu$ be a positive non-pluripolar measure on $X$. Then 
\begin{itemize}
\item[(i)] if in addtion $\mu(X)=\int_X \mbox{MA}(\phi)$,there exists $u\in\mathcal{E}(X,\beta,\phi)$ (unique up to a constant) such that $$(\beta+dd^c u)^n=\mu;$$
\item[(ii)] there exists $u\in\mathcal{E}^1(X,\beta,\phi)$ such that $$(\beta+dd^c u)^n=e^u\mu,$$
where $\mathcal{E}(X,\beta,\phi)$ is the relative full mass class defined in \ref{def of relative full mass class} and $\mathcal{E}^1(X,\beta,\phi)$ is the relative finite energy class defined in \ref{defn of finite energy class}.
\end{itemize}
\end{thm}

As an important application of the above theorem, we borrow Darvas-Di Nezza-Lu's idea in \cite{DDNL18a} and give an explicit description of Lelong numbers of elements in $\mathcal{E}(X,\beta,\phi)$.
\begin{thm}
Assume that $\varphi\in\mathcal{E}(X,\beta,\phi)$. Then $\forall x\in X$, 
$$\nu(\varphi,x)=\nu(\phi,x).$$
\end{thm}
\begin{rem}
If $\beta$ is semipositive and then $\phi=0$, this theorem is proved by Darvas-Di Nezza-Lu when $X$ is a compact K\"ahler manifold in \cite{DDNL18a}. Their result was a solution of an open problem of Dinew-Guedj-Zeriahi \cite[Problem 36]{DGZ16}: When $X$ is a compact K\"ahler manifold (or a manifold from the Fujiki class), $\beta$ is a semi-positive big form on $X$, is it true that any element of $\mathcal{E}(X,\beta)$ has zero Lelong number?

 In view of our theorem, the K\"ahler (or the Fujiki class) condition is shown to be unnecessary, which can be seen as a solution of a generalized version of Dinew-Guedj-Zeriahi's problem.
\end{rem}

\section{Pluripotential theory on compact Manifolds}
Let (X, $\omega$) be a compact Hermitian manifold of complex  dimension $n$, with a Hermitian metric $\omega$.
Let $\beta$ be a smooth real closed $(1,1)$ form. 
A  function $u:X\rightarrow [-\infty,+\infty)$ is called quasi-plurisubharmonic, if locally $u$ can be written as the sum of a smooth function and a psh function.  A  $\beta$-plurisubharmonic ($\beta$-psh for short) is a quasi-plurisubharmonic function $u$ such that    $\beta+dd^cu\geq 0$ in the sense of currents. 
		
The set of  all $\beta$-psh  functions on $X$ is denoted by $\mbox{PSH}(X,\beta)$. Suppose that  there exists a function  $\rho\leq 0 \in \mbox{PSH}(X,\beta)\cap L^{\infty}(X)$. Note that in this case the function $\sup\{\varphi|\varphi\in\mbox{PSH}(X,\beta)\leq0\}$ is also bounded and negative. For convenience  we denote it by $\rho$ if  not specifically stated. 

\subsection{Non-pluripolar product}
\begin{defn}\label{defn: nn product}
Let  $\varphi \in \mbox{PSH}(X,\beta)$, the non-pluripolar product $ \left\langle(\beta+dd^c \varphi)^n \right\rangle$ is defined as 
	$$(\beta+dd^c \varphi)^n :=\lim_{k\rightarrow \infty} \mathds{1}_{\{\varphi \textgreater \rho-k\}}(\beta+dd^c\varphi^{(k)})^n,$$
	where $\varphi^{(k)}:=\max\{\varphi,\rho-k\}$. 
\end{defn}

	It is easy to see that  $O_k:=\{\varphi>\rho-k\}$ is a plurifine open subset, and for any compact subset  $K$ of $X$, 
	\begin{align*}
		\sup_k\int_{K\cap O_k}(\beta+dd^c\varphi^{(k)})^n&\leq \sup_k\int_X(\beta+dd^c\varphi^{(k)})^n\\
		&=\int_X\beta^n<+\infty.
	\end{align*}
	Thus the Definition \ref{defn: nn product} of non-pluripolar product is well defined in the sense of \cite[Definition 1.1]{BEGZ10}.

\begin{defn} The non-pluripolar Monge–Ampère measure of $\varphi$ (with respect to $\beta$) is denoted by
 $MA(\varphi):=(\beta+dd^c \varphi)^n$.
    We say that  $\varphi$ has full Monge–Ampère mass if $\beta+dd^c \varphi$ has full Monge–Ampère mass, that is if and only if the measure MA($\varphi$) satisfies
    $$\int_{X}\mathrm{MA}(\varphi)=\int\beta^n.$$
\end{defn}
Using local plurifine property and following \cite{BEGZ10}, we can easily prove that
\begin{prop}\label{comp princi}
For two $\beta$-psh function $\varphi$ and $\psi$,we have
    $$\int_{\{\varphi<\psi\}}\mathrm{MA}(\psi)\leqslant\int_{\{\varphi<\psi\}}\mathrm{MA}(\varphi)+\int_X \beta^n-\int_{X}\mathrm{MA}(\varphi).$$
\end{prop}
\begin{defn}\label{BTcap}
Fixed $\mathcal{U}=\{\mathcal{U}_{\alpha}\}$ a finite covering of $X$ by strictly pseudoconvex open subsets of $X$,$\mathcal{U}_{\alpha}=\{x\in X:\rho_{\alpha}<0\}$, where $\rho_{\alpha}$ is a strictly psh smooth function defined in a neighborhood of $\overline{{{\mathcal{U}}_{\alpha}}}.$ Fix $\delta>0$ such that $\mathcal{U}^{\delta}=\{\mathcal{U}_{\alpha}^{\delta}\}$ is still a covering of $X$, where $\mathcal{U}_{\alpha}^{\delta}=\{x\in X:\rho_{\alpha}(x)<-\delta\}$. For a Borel subset $K$ of $X$, we define
$$Cap_{BT}(K):=\sum_{\alpha}Cap_{BT}(K\cap\mathcal{U}_{\alpha}^{\delta},\mathcal{U}_{\alpha}),$$
where $$Cap_{BT}(E,\Omega):=\sup\left\{\int_E(dd^cu)^n|u\in PSH(\Omega),0\leq u\leq1\right\}$$ is the capacity studied by Bedford and Taylor.

A function $u$ is called quasi-continuous if  for each $\varepsilon>0$, there exists an open set U such
that $\mbox{Cap}_{BT}(U)$ and the restriction of $u$ on $X-U$ is continuous.
A sequence of functions $u_j$ converges in capacity to $u$ if for any $\delta>0$,
$$\lim\limits_{j\to\infty}\text{Cap}_{BT}(\{x\in X:|u_j(x)-u(x)|>\delta\})=0.$$
\end{defn}
We can reformulate \cite[proposition 4.25]{GZ17} into a global version:
\begin{prop}
Let $(\varphi_j)_j\subset\mbox{PSH}(X,\beta)$ be a monotone sequence which converges almost everywhere to $\varphi\in\mbox{PSH}(X,\beta)$. Then $(\varphi_j)$ converges to $\varphi$ in capacity.
\end{prop}

Now we can prove a  lower-semicontinuity property of non-pluripolar products:
\begin{thm}\label{convergence thm}
Let $\beta^{j},j\in\{1,...,n\}$ be smooth closed real $(1,1)-$forms on $X$ whose cohomology classes have bounded potentials $\rho_j$. Suppose that for all $j\in\{1,...,n\}$ we have $u_j,u_j^k\in\mbox{PSH}(X,\beta^j)$ such that $u_j^k\to u_j$ in capacity as $k\to\infty$, and let $\chi_k,\chi\geq 0$ be quasi-continuous and uniformly bounded such that $\chi_k\to\chi$ in capacity. Then
\begin{equation}\label{1}
\liminf_{k\to\infty}\int_{X}\chi_{k}(\beta^1+dd^c u_1^k)\wedge ...\wedge (\beta^n+dd^c u_n^k)\geq\int_X \chi(\beta+dd^c u_1)\wedge...\wedge (\beta^n+dd^c u_n).
\end{equation}
If in addtion,
\begin{equation}\label{2}
\int_X (\beta^1+dd^c u_1)\wedge...\wedge (\beta^n+dd^c u_n)\geq \limsup_{k\to\infty}\int_X (\beta^1+dd^c u_1^k)\wedge...\wedge (\beta+dd^c u_n^k),
\end{equation}
then $\chi_k (\beta^1+dd^c u_1^k)\wedge...\wedge (\beta^n+dd^c u_n^k)\to\chi (\beta^1+dd^c u_1)\wedge...\wedge (\beta+dd^c u_n)$ in the weak sense of measures on $X$.
\end{thm}
\begin{proof}
Fix an open relatively compact subset $U$ of $X$ and replace $V_{\theta}$ in the proof of \cite[Theorem 2.6]{DDNL23} by $\rho_j$. Follow the lines and we will obtain the proof.
\end{proof}
\begin{rem}\label{increase implies 2}
If $u_{j}^{k}\nearrow u_{j}$ a.e. as $k\to\infty$, then $(2)$ is automatically satisfied. See \cite[Remark 3.4]{DDNL23}.
\end{rem}
We also define the Monge-Amp\`ere capacity with respect to $\beta$ as below and prove that a subset is (locally) pluripolar if and only if it is $\mbox{PSH}(X,\beta)$-polar, which can be seen as a generalization of \cite{GZ05}.  For every Borel subset $B$ of $X$,
\begin{equation}
     \mathrm{Cap}_{\beta}(B):=\sup\biggl\{\int_{\mathrm{B}}\mathrm{MA}(\varphi)\mid\varphi\in\mathrm{PSH}(\mathrm{X},\beta),\rho-1\le\varphi\le\rho ~\mathrm{on}~\mathrm{X}\biggr\}
 \end{equation}
See \cite{GZ05} for the K\"ahler case and \cite{BBGZ13} for the big case.
\begin{prop}\label{basic properties of cap}
The capacity defined above enjoys basic properties as follows:
\begin{itemize}
\item[(i)]  If $K$ is a compact subset, the supremum in the definition of $\mathrm{Cap(K)}$ is achieved by the usc regularization of
    \begin{align*}
        h_{\mathrm{K}}:&=\sup\bigl\{\varphi\in\mathrm{PSH}(\mathrm{X},\beta)\mid\varphi\leq\rho ~on~\mathrm{X},\varphi\leq\rho-1~ on~\mathrm{K}\bigr\}\\
   &=\sup\bigl\{\varphi\in\mathrm{PSH}(\mathrm{X},\beta)\mid\varphi\leq 0~ on~\mathrm{X},\varphi\leq\rho-1~ on~\mathrm{K}\bigr\}.
   \end{align*}

\item[(ii)] Let $E$ be a Borel subset. Then we have
$h_{\mathrm{E}}^{*}=\rho-1$ a.e. on $E$ and $h_{\mathrm{E}}^{*}=\rho$ a.e.on $X-\Bar{E}$ with respect to the measure $\mathrm{MA}(h_{\mathrm{E}}^{*})$.\\
\item[(iii)] Let $K$ be a compact set of $X$, $\{K_j\}$ be a sequence of compact sets decrease to $K$. Then 
$$\mathrm{Cap}_{\beta}(K)=\lim_{j\rightarrow\infty}\mathrm{Cap}_{\beta}(K_j).$$
In particular, $$\mathrm{Cap}_{\beta}(K)=\inf\{\mathrm{Cap}_{\beta}(U):K\subseteq U, \ U \ \mathrm{is} \ \mathrm{open}\}.$$
\end{itemize}
\end{prop}
\begin{proof}
For $(i)$, the proof is similar to  \cite[Lemma 1.4]{BBGZ13}. It is easy to see that $h_K^{*}$ is a candidate in the definition. Pick any arbitrary $\varphi\in\mbox{PSH}(X,\beta)$ satifying $\rho-1\leq\varphi\leq\rho$. It is enough to show that $\int_K (\beta+dd^c \varphi)^n\leq\int_K (\beta+dd^c h_K^{*})$. By definition we have $h_K^{*}<(1-\epsilon)\varphi+\epsilon\rho$ and $K\subset\{h_K^{*}<(1-\epsilon)\varphi+\epsilon\rho\}$. Then
\begin{align*}
\int_K (\beta+dd^c (1-\epsilon)\varphi+\epsilon\rho)^n\leq\int_{\{h_K^{*}<(1-\epsilon)\varphi+\epsilon\rho+1\}}(\beta+dd^c \varphi)^n&\leq\int_{\{h_K^{*}<(1-\epsilon)\varphi+\epsilon\rho\}}(\beta+dd^c h_K^{*})^n\\
&\leq\int_{\{h_K^{*}<\rho\}}(\beta+dd^c h_K^{*})^n=\int_K (\beta+dd^c h_K^{*})^n.
\end{align*}
The third inequality holds by comparison principle \ref{comp princi} (bounded $\beta$-psh function always has full Monge-Amp\`ere mass). The final equality uses $(ii)$.
By mulitilinearity of non-pluripolar product we let $\epsilon\to0$ and we obtain the result.

 For $(ii)$, note that $h_{E}\leq\rho-1\leq h_{E}^{*}~\mbox{on~K}$. By Bedford-Taylor's theorem the measure $\mbox{MA}(h_E^{*})$ does not charge the pluripolar set $\{h_E<h_E^{*}\}$, and the first claim holds. To prove the second, one need to prove $\mathrm{MA}(h_{\mathrm{E}}^{*})$ puts no mass on $\{h^{*}_K < \rho\}-\Bar{E}\subset\{h^{*}_K < 0\}-\Bar{E}$. Similarly to the proof of \cite[Proposition 4.1]{GZ05}, $\mathrm{MA}(h_{\mathrm{E}}^{*})$ puts no mass on $\{h^{*}_K < 0\}-\Bar{E}$.

 For $(iii)$, the proof is similar to \cite[Corollary 4.6]{DDNL18b}. For the reader's convenience we repeat it here. Set $h_j:=h_{K_j}^*$. It is obvious that $\{h_j\}$ increases to a $\beta$-psh function $h$ a.e. on $X$. We are going to prove that $h=\rho-1$ a.e. on $K$ and $h=\rho=0$ a.e. on $X-K$ with respect to $\mathrm{MA}(h)$. The first claim is obvious since $h=\rho-1$ on $K-\cup\{h_{K_j}\textless h_{K_j}^*\}$ and $\{h_{K_j}\textless h_{K_j}^*\}$ is pluripolar. To prove the second claim, note that $\{h\textless0\}-K_j\subseteq \{h_m\textless0\}-K_m$ when $m\geq j$, we get for fixed $j$,
$$\mathrm{MA}(h)(\{h\textless0\}-K_j)\leq \liminf_{m\rightarrow\infty}\mathrm{MA}(h_m)(\{h\textless0\}-K_j)\leq \mathrm{MA}(h_m)(\{h_m\textless0\}-K_m\})=0.$$
 Let $j\rightarrow\infty$ and the second claim follows. It then follows from $(i)$ and Remark \ref{increase implies 2} that
$$\lim_{j\rightarrow\infty}\mathrm{Cap}_{\beta}(K_j)=\lim_{j\rightarrow\infty}\int_X(h_j-\rho)\mathrm{MA}(h_j)=\int_X(h-\rho)\mathrm{MA}(h)=\int_K\mathrm{MA}(h)\leq \mathrm{Cap}_{\beta}(K).$$
Since the reverse inequality is obvious, the first statement follows. To prove the second statement, we take $\{K_j\}$ to be a decreasing sequence of compact neighbourhood of $K$. By the first claim we get that
$$\mathrm{Cap}_{\beta}(K)=\lim_{j\rightarrow\infty}\mathrm{Cap}_{\beta}(K_j)\geq \lim_{j\rightarrow\infty}\mathrm{Cap}_{\beta}(\mathring{K_j})\geq \inf\{\mathrm{Cap}_{\beta}(U):K\subseteq U, \ U \ \mathrm{is} \ \mathrm{open}\}.$$ The reverse inequality is obvious and the second statement follows.

\end{proof}
\noindent{\bf Two capacities:} We now introduce two capacities as what has been done in \cite{GZ05} and prove a more general theorem which reveals the equivalence between (locally) pluripolarity and $\mbox{PSH}(X,\beta)$-polarity. This is important in our proof of solving degenerate Monge-Amp\`ere equations with respect to a cohomology class having bounded potentials on compact Hermitian manifolds.

\begin{defn}
Let $E$ be a Borel set of $X$ and $\rho^{\prime}:=\rho-\inf_X\rho$, we define the global $\beta$-extremal functions  by
$$V_{K,0}=\sup \{\varphi\in \mbox{PSH}(X,\beta):\varphi\leq0 \ \mathrm{on} \ K \},$$
$$V_{K,\rho}=\sup \{\varphi\in \mbox{PSH}(X,\beta):\varphi\leq \rho^{\prime} \ \mathrm{on} \ K \}.$$
We then introduce the Alexander-Taylor capacities of $E$,
$$T_{\beta,0}(E)=\exp(-\sup_{X}V_{K,0}^*),$$
$$T_{\beta,\rho}(E)=\exp(-\sup_{X}V_{K,\rho^{\prime}}^*).$$
\end{defn}

We are going to prove that every (local) pluripolar set is actually $\mbox{PSH}(X,\beta)-$polar. To do this, we need the following proposition.
\begin{prop}\label{prop:compact set exp control}
There exists $C\textgreater0$ such that for any compact subset $K$ of $X$,
$$T_{\beta,\rho}(K)\geq\exp(-\frac{C}{\mbox{Cap}_{\beta}(K)}); T_{\beta,0}(K),T_{\beta,\rho^{\prime}}(K)\leq C\exp(-\frac{\emph{Cap}_{\beta}(K)^{-\frac{1}{n}}}{C}).$$
\end{prop}
\begin{proof}
We prove the latter two inequality first. The proof of them are the same, so we simply write once.If $\sup_XV_{K}^*=+\infty$, then $K$ is pluripolar, $T_{\beta}(K)=\exp(-{\mathrm{Cap}_{\beta}(K)^{-\frac{1}{n}}})=0$, thus there is nothing to prove. So we assume $M_K:=\sup_XV_{K}^*\textless +\infty$ and set $M=||\rho^{\prime}||_{\infty}$. There are two cases to deal with: \\
Case 1: $M+M_K\geq1$, we get that
\begin{align*}
\frac{1}{(M_K+M)^n}\int_X\beta^n&=\frac{1}{(M_K+M)^n}\int_X(\beta+dd^cV_K^*)^n\\&=\frac{1}{(M_K+M)^n}\int_K(\beta+dd^cV_K^*)^n\\&\leq\int_K(\beta+dd^c(\frac{V_K^*}{M_K+M}+\frac{(M_K+M-1)\rho^{\prime}}{M_K+M}))^n.
\end{align*}
By the definition of $\mathrm{Cap}_{\beta}(\cdot)$ and the fact that
$$\rho^{\prime}-\frac{M}{M_K+M}\leq\frac{V_K^*}{M_K+M}+\frac{(M_K+M-1)\rho^{\prime}}{M_K+M}\leq \rho^{\prime}+\frac{M_K}{M_K+M},$$
 we get that $(M_K+M)^n\int_X\beta^n\leq Cap_{\beta}(K)$, which is equivalent to $$T_{\beta}(K)\leq \exp(-M-(\frac{\mathrm{Cap}_{\beta}(K)}{\mathrm{vol}(\beta)})^{-\frac{1}{n}}).$$
Case 2: $0\leq M_K+M\leq1$, we get that
\begin{align*}
\int_X\beta^n
&=\int_X(\beta+dd^cV_K^*)^n\\
&=\int_K(\beta+dd^cV_K^*)^n.
\end{align*}
By the definition of $\mathrm{Cap}_{\beta}(\cdot)$ and the fact that $\rho^{\prime}-M\leq V_K^*\leq\rho^{\prime}+M_K$, we get that $\int_X\beta^n\leq \mathrm{Cap}_{\beta}(K)$, this also leads to
$$T_{\beta}(K)\leq \exp(-M-(\frac{\mathrm{Cap}_{\beta}(K)}{\mathrm{vol}(\beta)})^{-\frac{1}{n}}).$$
For the first inequality,  let $\varphi\in\mbox{PSH}(X,\beta)$ be such that $\varphi\leq \rho$ on $K$. Then $\frac{\varphi-M_K-M}{M_K+M}+\frac{(M_K+M-1)\rho}{M_K+M}\in\mbox{PSH}(X,\beta)$ satisfies $\leq \rho-1$ on $K$ and $\leq 0$ on X. Thus
$$\frac{V_K^{*}-M_K-M}{M_K+M}+\frac{(M_K+M-1)\rho}{M_K+M}\leq h_K^{*}\leq0.$$
Now $\sup_X({V_K^{*}-M_K})=0$, it follows that there exists a constant $C$ independent of K such that $\int_X |V_K^{*}-M_K|\omega^n \leq C$. We infer that
\begin{align*}
\mbox{Cap}_{\beta}(K)=\int_K (\rho-h_{K}^{*})\mbox{MA}(h_{K}^{*})&\leq\int_X (\rho-h_K^{*})\mbox{MA}(h_K^{*})\\&\leq \frac{1}{M_K+M}\int_X ({\rho+M+M_K-V_K^{*}})\mbox{MA}(h_K^{*})\leq \frac{C}{M_K}
\end{align*}
using Proposition \ref{basic properties of cap} $(i)$ and CLN-inequality for bounded $\beta$-psh functions. This yields the desired inequality.
\end{proof}
\begin{rem}
When $\beta$ is K\"ahler, it is obvious that $V_{K,0}=V_{K,\rho^{\prime}}=V_{K,\rho}$ and $M=0$. Thus Proposition \ref{prop:compact set exp control} can be seen as a generalization of \cite[Proposition 7.1]{GZ05}
\end{rem}
\begin{cor}\label{coro:open set exp control}
There exists $C\textgreater0$ such that for any open subset $G$ of $X$,
$$T_{\beta,\rho^{\prime}}(G)\leq C\exp(-\frac{\mathrm{Cap}_{\beta}(G)^{-\frac{1}{n}}}{C}).$$
\end{cor}
\begin{proof}
Let $K_j$ be increasing sequence of compact sets such that $K_j$ and $\mathring{K_j}$ exhaust $G$. Since $V_{K_j,\rho^{\prime}}=\rho^{\prime}$ on $\mathring{K_j}$, we know that the limit $V:=\lim_{j\rightarrow \infty}V_{K_j,\rho^{\prime}}^*=\rho^{\prime}$ on $G$. Hence $V$ is a candidate of $V_{G,\rho^{\prime}}$, we get that $V\leq V_{G,\rho^{\prime}}$. On the other hand, by definition, $V_{K_j,\rho^{\prime}}^*\geq V_{G,\rho^{\prime}}^*$, we get that $V\geq V_{G,\rho^{\prime}}^*$. This leads to $V=V_{G,\rho^{\prime}}=V_{G,\rho^{\prime}}^*$.
\end{proof}
Now we can prove the following theorem stated at the begining.
\begin{thm}\label{thm:local is global}
Every $(locally)$ pluripolar set is $\mbox{PSH}(X,\beta)$-polar.
\end{thm}
\begin{proof}
The proof is similar to \cite[Theorem 7.2]{GZ05}. One only need to replace their estimate by the one in Corollary \ref{coro:open set exp control}.
\end{proof}

\begin{thm}

$(1)$ Let $E$ be a Borel set and $P$ be a $\mbox{PSH}(X,\beta)$-polar set, then
$V_{E,0}^*=V_{E-P,0}^*.$\\
$(2)$ Let $\{ E_j\}$ be a increasing sequence of Borel sets and set $E=\bigcup E_j$.Then $V_{E,0}^*=\lim_{j\rightarrow \infty}V_{E_j,0}^*.$\\
$(3)$ Let $\{ K_j\}$ be a decreasing sequence of compact sets and set $K=\bigcap K_j$.Then $V_{K_j,0}^*$ increases a.e. to $V_{K,0}^*$.
\end{thm}

\begin{proof}

\noindent{\bf Proof of assertion (1):} It is obvious that $V_{E,0}^*\leq V_{E-P,0}^*$. To prove the reverse, assume $\psi\in \mbox{PSH}(X,\beta)$ such that $\psi\leq0$ and $N\subseteq\{\psi=-\infty\}$. Fix a candidate of $V_{E-P,0}$, namely $\varphi$. Then $(1-\varepsilon)\varphi+\varepsilon \psi\leq 0$ on $E$ for any $0 \textless \varepsilon\textless 1$. Therefore $(1-\varepsilon)\varphi+\varepsilon \psi\leq V_{E,0}^*$. Let $\varepsilon\rightarrow0$ we get that $\varphi\leq V_{E,0}^*$ on $X-\{\psi=-\infty\}$, hence on the whole $X$. Therefore $V_{E-P,0}^*\leq V_{E,0}^*$.\\
\noindent{\bf Proof of assertion (2): }Since $E_j$ is increasing to $E$, one can see $\{V_{E_j,0}^*\}$ is decreasing and $V_{E_j,0}^*\geq V_E^*$. Let $V=\lim_{j\rightarrow \infty}V_{E_j,0}^*$. If $E$ is $\mbox{PSH}(X,\beta)$-polar, then $V_{E,0}^*=V_{E_j,0}^*=+\infty$, so we may assume $E$ is not $\mbox{PSH}(X,\beta)$-polar. In this case, $V$ and $V_{E,0}^*$ are both $\beta$-psh functions. Observe that $V\leq0$ on $E-N$, where $N=\cup_{j=1}^{\infty} \{V_{E_j,0}\textless V_{E_j,0}^* \}$ is a (locally) pluripolar set by \cite{BT82}. By Theorem \ref{thm:local is global} above we get that $N$ is a $\mbox{PSH}(X,\beta)$-polar set. Therefore $V_{E,0}^*\leq V\leq V_{E-N,0}^*=V_{E,0}^*$ by assertion (1).\\
\noindent{\bf Proof of assertion (3):} Take a candidate of $V_{K,0}$, namely $\varphi\in \mbox{PSH}(X,\beta)$ such that $\varphi \leq 0$ on $K$. Then $\{\varphi\textless \varepsilon\}$ is an neibourhood of $K$ which will contains $K_j$ when $j\textgreater j_{\varepsilon}$. Hence $\varphi-\varepsilon\leq \lim_{j\rightarrow\infty}V_{K_j,0}$. Since $\varphi$ and $\varepsilon$ are arbitrary, we get that $V_{K,0}\leq \lim_{j\rightarrow \infty}V_{K_j,0}$. On the other hand, $V_{K_j,0}\leq V_{K,0}$ is obvious. Thus we get the equality, the desired result follows.
\end{proof}

With assetions (2),(3) above, one can easily see that the set funtion $T_{\beta,0}$ is a capacity:
\begin{prop}
\begin{itemize}
\item[(i)]The set function $T_{\beta,0}(\cdot)$ satisfies:
$$T_{\beta,0}(E_1)\leq T_{\beta,0}(E_2)~for~ E_1\subseteq E_2;$$
\item[(ii)]If $\{K_j\}$ is an decreasing sequence of compact sets and $K=\bigcap_{j=1}^{\infty} K_j$, then $$T_{\beta,0}(K)=\lim_{j\rightarrow \infty}T_{\beta,0}(K_j);$$

\item[(iii)]If $\{ E_j\}$ be a increasing sequence of Borel sets and $E=\bigcup_{j=1}^{\infty} E_j$, then $$T_{\beta,0}(E)=\lim_{j\rightarrow \infty}T_{\beta,0}(E_j).$$

\end{itemize}
\end{prop}

The well-known  Choquet capacitability theorem  implies following theorem:
\begin{thm}
Every Borel subset $E$ of $X$ satisfies
$$T_{\beta,0}(E)=\sup_{K \subseteq \subseteq E}T_{\beta,0}(K).$$
\end{thm}
As a corollary, we get that
\begin{cor}\label{coro:strong exp control}
There exists $C\textgreater0$ such that for any Borel subset $E$ of $X$,
$$T_{\beta,0}(E)\leq C\exp(-\frac{\mathrm{Cap}_{\beta}(E)^{-\frac{1}{n}}}{C}).$$
\end{cor}

\begin{thm}\label{local polar=global polar}
Let $E$ be a Borel set of $X$, then the following statement are equivalent:\\
$(1)$ E is $\mbox{PSH}(X,\beta)$-polar;
$(2)$ E is $(local)$ pluripolar;
$(3)$ $\mathrm{Cap}_{BT}(E)=0$;
$(4)$ $\mathrm{Cap}_{\beta}(E)=0$.
\end{thm}
\begin{proof}
$(1)\Rightarrow(2)\Rightarrow(3)\Rightarrow(4)$ is obvious. If $(4)$ holds, i.e. $Cap_{\beta}(E)=0$, we get that $\sup_XV_{E,0}^*=+\infty$ by Corollary \ref{coro:strong exp control}. Then it follows from almost same argument in \cite[Theorem 5.2]{GZ05} that $(1)$ holds.
\end{proof}

\subsection{Monotonicity of Non-pluripolar product masses}
Given two $\beta$-psh functions $\phi$ and $\psi$, we say that $\phi$ is more singular than $\psi$ if there exists a constant $C>0$ such that $\phi\leq\psi+C$. We say that they have the same singularity type if there exists a constant $C>0$ such that $-C\leq \phi-\psi\leq C$. In this part we will prove that $\beta$-psh function with more singularity has less Monge-Amp\`ere mass. The ideas are adapted from \cite{DDNL23} where  quasi-psh potentials from a big cohomology class  are dealed. See also \cite{WN19}\cite{Vu21}.
\begin{lem}\label{same sing same mass}
Let $u,v\in\mbox{PSH}(X,\beta)$. If $u$ and $v$ have the same singularity type, then $\int_X (\beta+dd^c u)^n=\int_X(\beta+dd^cv)^n$.
\end{lem}
\begin{proof}
The proof is similar to \cite[Lemma 3.1]{DDNL23} with some modification and is divided into two steps. We may assume that $\rho>0$. Firstly by the proof of \cite{DDNL23} we have that: if there exists a constant $C>0$ such that $u=v$ on open set $U:=\{\min(u,v)<C\}$ then $\int_X (\beta+dd^c u)^n=\int_X(\beta+dd^cv)^n$.

Now since $u$ and $v$ have same singularity type, there exists a positive constant $B$ such that $v\leq u\leq v+B\leq 0$. For each $a\in (0,1)$ we set $v_a:=av+(1-a)\rho,u_a:=\max (u,v_a)$ and $ C=Ba(1-a)^{-1}>Ba(1-a)^{-1}-(1-a)\rho$. It is easy to check that $u_a=v_a$ on the open set $U_t:=\{\min (u_a,v_a)<-C\}$. Applying the claim we get
 $\int_X \mbox{MA}(u_a)=\int_X \mbox{MA}(v_a)$. Since non-pluripolar products are multilinear, we have that
$$\int_X \mbox{MA}(v_a)=a^n\int_X \mbox{MA}(v)+\sum_{k=0}^{n-1}a^k(1-a)^{n-k}\int_X (\beta+dd^c v)^{k}\wedge (\beta+dd^c \rho)^{n-k}\to\int_X \mbox{MA}(v)$$
as $a$ increase to 1. Since $u_a$ decrease to $u$ as $a$ increase to 1, by Theorem \ref{convergence thm} we have
$$\liminf_{a\to 1^{-}}\int_X \mbox{MA}(u_a)\geq\int_X \mbox{MA}(u).$$
We thus have $\int_X \mbox{MA}(u)\leq\int_X \mbox{MA}(v).$ By symmetry we get equality and finish the proof.
\end{proof}
\begin{rem}
Although we are dealing with quasi-psh functions in the class with bounded potentials, it is easy to see that this also implies \cite[Lemma 3.1]{DDNL23}. Indeed, when $\beta$ is K\"ahler and $\theta$ is a smooth representative $(1,1)-$form of a big class on $X$, we let $t$ big enough such that $\theta+t\beta$ is K\"ahler. If $u,v\in\mbox{PSH}(X,\theta)$ and have the same singularity type,then by Lemma \ref{same sing same mass} $\int_X (\theta+t\beta+dd^c u)^n=\int_X (\theta+t\beta+dd^c v)^n$ for $t$ big enough hence for all $t>0$ by multi-linearity of the non-pluripolar product. Comparing the coefficience, we have that $\int_X (\theta+dd^c u)=\int_X (\theta+dd^c v)$.

\end{rem}
Using the above result and adapt the proof of \cite{DDNL23}, we  can also obtain that
\begin{prop}\label{multi mass eq}
Let $\beta^j,j\in\{1,...,n\}$ be smooth closed real $(1,1)-$forms on $X$ whose cohomology classes have bounded potentials $\rho_j$. Let $u_j,v_j\in\mbox{PSH}(X,\beta^j)$ such that $u_j$ has the same singularity type as $v_j,j\in\{1,...,n\}$. Then
$$\int_X (\beta^1+dd^c u_1)\wedge...\wedge (\beta^n+dd^c u_n)=\int_X (\beta^1+dd^c v_1)\wedge...\wedge (\beta^n+dd^c v_n).$$
\end{prop}
\begin{thm}\label{mass ineq}
Let $\beta^j,j\in\{1,...,n\}$ be smooth closed real $(1,1)-$forms on $X$ whose cohomology classes have bounded potentials $\rho_j$.Let $u_j,v_j\in\mbox{PSH}(X,\beta^j)$ be such that $u_j$ is less singular than $v_j$ for all $j\in\{1,...,n\}$. Then
$$\int_X (\beta^1+dd^c u_1)\wedge...\wedge (\beta^n+dd^c u_n)\geq\int_X (\beta^1+dd^c v_1)\wedge...\wedge (\beta^n+dd^c v_n).$$
\end{thm}

\subsection{Envelopes}
\begin{defn}
If $f$ is a function on $X$, we define the envelope of $f$ in $\mbox{PSH}(X,\beta)$ by $$P(f):=(\sup\{u\in\mbox{PSH}(X,\beta)|u\leq f\})^{*}$$ with the convention that $\sup\emptyset=-\infty$.
\end{defn}
It is easy to see that if $f$ is usc, $P(f)$ is finite and satisfies $P(f+C)=P(f)+C$ for any constant $C$. If $f=\min(\psi,\phi)$ where $\psi,\phi\in\mbox{PSH}(X,\beta)$, we denote $P(min(\psi,\phi))$ by $P(\psi,\phi)$. 
\begin{prop}\label{mass concentr of envelope}
    Let $u:X\to\mathbf{R}$ be a quasi-continuous and usc function. Then $\int_{\{P(u)<u\}}(\beta+dd^c P(u))^n =0$.
\end{prop}
\begin{proof}
     First of all we treat the case that $u$ is continuous. By Choquet's lemma, there exist $(u_j)\subset\mbox{PSH}(X,\beta)\cap L^{\infty}(X)$ such that $(\lim u_j)^{*}=P(u)$. Let $x_0\in\{P(u)<u\}$. Since $u$ is continuous, there exists a small open ball $B=B(x_0,r)$ such that for $x\in \overline{B(x_0,r)}$, $P(u)(x)<u(x)-2\varepsilon<u(x_0)-\varepsilon<u(x)$. By classical pluripotential theory, we can find $\varphi_j\in\mbox{PSH}(X,\beta)$ such that $\varphi_j=u_j$ on $X \setminus B$, $\mbox{MA}(\varphi_j)=0$ in $B$, and $\varphi_j$ is still increasing. 

To treat the general case, we borrow the idea from \cite[Proposition 2.16]{DDNL18a}. By semi-continuity we can approximate $u$ from above by a sequence of smooth functions $(f_j)$. Set $\varphi_j:=P(f_j),\varphi:=P(u)$ and note that $\varphi_{j}\searrow\varphi.$ For each $j\in\mathbb{N}$ the measure $\mbox{MA}(\varphi_j)$ vanishes in the set $\{\varphi_j<f_j\}$. Now we want to pass to the limit as $j\to\infty$. Fix $k,l\in\mathbb{N}$ and set $$U_{k,l}=\{\varphi_k<u\}\cap\{\varphi>\rho-l\}.$$
For any $j>k$,note that on $U_{k,l}$ we have $\varphi_j\geq\varphi>\rho-l$ and $\{\varphi_k< u\}\subset\{\varphi_j<f_j\}$. By plurifine property of non-pluripolar product, for any $j>k$, the measure $\mbox{MA}(\max(\varphi_j,\rho-l))=\mbox{MA} (\varphi_j)$ vanishes on $U_{k,l}$. By assumption $u$ is quasi
continuous (see Definition \ref{BTcap}), hence $U_{k,l}$ is quasi open. More precisely, for any fixed $\epsilon>0$ there exists an
 open set $V_\epsilon$ such that the set $G_{\varepsilon}:=(V_{\varepsilon}\setminus U_{k,l}\cup U_{k,l}\setminus V_{\varepsilon})$ satisfies $\mathrm{Cap}_{BT}(G_{\varepsilon})\leq\varepsilon $ (see Definition \ref{BTcap}). Observe that for fixed $l$ ,$\rho-l\leq\max(\varphi_j,\rho-l)\leq\rho$. It follows that 
$$\sup_{j\in\mathbb{N}}\int_{G_{\epsilon}} \operatorname{MA}({\max(\varphi_{j},\rho-1)})\leq A\text{Cap}_{\beta}(G_{\epsilon})\leq A'\text{Cap}_{BT}(G_{\epsilon})$$ where the last inequality follows from  definitions of two capacities. \\
Consequently, $\sup_{j\in\mathbb{N}}\int_{V_{\epsilon}}\operatorname{MA}({\max(\varphi_{j},\rho-l)})\leq A^{\prime}\varepsilon$. Take the limit and we ultimately obtained:
$$\int_{U_{k,l}}\operatorname{MA}({\max(\varphi,\rho-l)})\leq C\varepsilon,$$ for some positive constant $C$ independent of $\epsilon$ (but dependent on $l$). Now letting $\epsilon\to 0$ we see that $\mbox{MA}\max(\varphi, \rho-j)$ vanishes in $U_{k,l}$. Letting $l\to +\infty$, and by definition of the non-pluripolar product, we see that $\mbox{MA}(\varphi)$ vanishes in $\{\varphi_k<u\}$. Now, letting $k\to +\infty$ we obtain the result.
\end{proof}
\begin{rem}
If $u$ is quasi-continuous, lower semi-continuous and satisfies $P(u)\in\mbox{PSH}(X,\beta)$, it is easy to deduce that $\mbox{MA}(P(u))$ is also supported on $\{P(u)=u\}$. Indeed, we can approximate $u$ from below by a sequence of smooth functions $(f_j)$. Let $\varphi_j:=P(f_j)$, then we have $\varphi_j$ increases to $\varphi$ up to a pluripolar subset (of BT-capacity zero). To deal with this case, it is enough to minus this subset at the suitable place and follow the line above.
\end{rem}
 As a byproduct of above remark, we can formulate \cite[Theorem 2.7]{DDNL23} as follows:
\begin{prop}
Let $u:X\to\mathbf{R}$ be a quasi-continuous and $P(u)\in\mbox{PSH}(X,\beta)$. Then $\int_{\{P(u)<u\}}(\beta+dd^c P(u))^n =0$.
\end{prop}

\noindent{\bf Regularity of quasi-psh envelopes:} When $\beta$ is K\"ahler and $\theta$ is a smooth representative $(1,1)-$form of a big class, Berman and Demailly proved the envelope $V_{\theta}:=\sup{\{u\in\mbox{PSH}(X,\theta),u\leq0\}}$ has locally bounded Laplacian in the Ample locus $\mbox{Amp}(\{\theta\})$ and its complex Monge-Ampère measures satifies $\mbox{MA}(V_{\theta})=\mathbf{1}_{\{V_{\theta}=0\}}\theta^n$. In the case when $\{\theta\}$ is integral, this result has been obtained by Berman \cite{Ber08} using different methods.We  borrow the idea of \cite[Theorem 2.6]{DDNL18a} and prove the regularity theorem below:
\begin{thm}\label{regularity envelop}
Let $\beta$ be a smooth $(1,1)-$form such that $\{\beta\}$ has a bounded potential. Then the envelope $\rho$ satisfies 
$$\mbox{MA}(\rho)\leq \mathbf{1}_{\{\rho=0\}}\beta^n.$$
\end{thm}
\begin{proof}
Set $\beta_{+}^n:=\max(\frac{\beta^n}{\omega^n},0)\omega^n$. By the result in \cite{LWZ23}, for each $\epsilon>0,\alpha>1$ there exists a unique bounded $\varphi_{\alpha,\epsilon}$ such that
$$\mbox{MA}(\varphi_{\alpha,\epsilon})=e^{\beta\varphi_{\alpha,\epsilon}}[(1+\epsilon)\beta_{+}^n+\epsilon\omega^n].$$ We claim that $\varphi_{\alpha,\epsilon}\leq0$ for all $\alpha>1,\epsilon>0$. By contradiction, we assume that there exists $x_0\in X$ such that $c:=\sup \varphi_{\alpha,\epsilon}>0$. Fix a ball $B(x_0,r)\subset X$ in a holomorphic coordinate chart around $x_0$. By shrinking $B(x_0,r)$, we can assume that the local potential $g$ of $\beta$ satifies $-c<g\leq0$ in $B(x_0,r)$. In this ball the function $u:=\varphi_{\beta,\epsilon}+g$ is plurisubharmonic, bounded, and we can write the following sequence of estimates:
\begin{align*}
(dd^cu)^n&=(\beta+dd^c \varphi_{\alpha,\epsilon})^n=e^{\alpha\varphi_{\alpha,\epsilon}}[(1+\epsilon)\beta^n_{+}+\epsilon\omega^n]\\
&\geq e^{\alpha u}[\beta^n_{+}+\epsilon\omega^n],
\end{align*}
where in the last line we used that $g\leq 0$ in $B(x_0,r)$ and $1+\epsilon>1$.
On the other hand the function $u-g-c$, defined in $B(x_0,r)$, attains a maximum at $x_0$ (equal to zero). It follows from \cite[Lemma 2.8]{DDNL18a} that at $x_0$ we have
$$\beta^n=(dd^c g)^n\geq e^{g+c}[\beta^n_{+}+\epsilon\omega^n].$$ But this is a contradiction since $g+c\geq 0$ and $\omega^n>0$ in $B(x_0,r)$. Thus $\varphi_{\alpha,\epsilon}\leq0$.
By \cite[Lemma 2.5]{DDNL18a}(reformulated in our setting) the solution $\varphi_{\alpha,\epsilon}\leq0$ are increasing as $\varepsilon\searrow0$. Indeed, if $s<\epsilon$ then
$$\mbox{MA}(\varphi_{\alpha,\epsilon})\geq e^{\alpha\varphi_{\alpha,\epsilon}}[(1+s)\beta^n_{+}+s\omega^n],$$
hence we can use \cite[Lemma 2.5]{DDNL18a}(reformulated in our setting) with $\mu:=(1+s)\beta^n_{+}+s\omega^n$ and $\phi=0$ to conclude that $\varphi_{\alpha,\epsilon}\leq\varphi_{\alpha,s}$.
Since $\varphi_{\alpha,\epsilon}\leq0$ for all $\alpha,\epsilon>0$ it follows that $\varphi_{\alpha,\epsilon}$ increase almost everywhere to some $0\geq\varphi_{\alpha}\in\mbox{PSH}(X,\beta)$ which is bounded on $X$ such that
$$\mbox{MA}(\varphi_{\alpha})=e^{\alpha\varphi_{\alpha}}\beta^n_{+}.$$
By \cite{LWZ23} again, there exists a unique  bounded $\phi\in\mbox{PSH}(X,\beta)$ such that
$$\mbox{MA}(\phi)=e^{\phi}\beta^n_{+}.$$
By \cite[Lemma 2.5]{DDNL18a}(reformulated in our setting) for $\alpha>1$ we have that
$$\varphi_{\alpha}\geq u_{\alpha}:=(1-\frac{1}{\alpha})\rho+\frac{1}{\alpha}\phi-\frac{n\log\alpha}{\alpha}.$$
Indeed, $u_{\alpha}\in\mbox{PSH}(X,\beta)$ is bounded and 
$$\mbox{MA}(u_{\alpha})\geq\frac{1}{\alpha^n}\mbox{MA}(\phi)=\frac{1}{\alpha^n}e^{\phi}\beta^n_{+}=e^{\phi-n\log\alpha}\beta^n_{+}\geq e^{\alpha u_{\alpha}}\beta^n_{+}.$$
It follows from \cite[Lemma 2.5]{DDNL18a}(reformulated in our setting) that $u_{\alpha}\geq\varphi_{\alpha}$ as claimed.
Since $\varphi_{\alpha}\geq0$, it follows  from the comparison principle that $\varphi_{\alpha}$ is increasing in $\alpha$. Hence $\varphi_{\alpha}\nearrow\rho$ and by continuity of the Monge-Ampère operator we have
$$\mbox{MA}(\rho)=\lim_{\alpha\to {+\infty}}(\beta+dd^c \varphi_{\alpha})^n\leq\beta^n_{+}.$$
Finally, we can prove that $\mbox{MA}(\rho)$ is supported on the contact set $\{\rho=0\}$. Indeed, for each $\delta>0$, since $U:=\{\rho<-\delta\}$ is open, we have
$$\int_{U}(\beta+dd^c \rho)\geq\liminf_{\alpha\to +\infty}\int_U e^{-\alpha\delta}\beta^n_{+}=0.$$
Hence $\mbox{MA}(\rho)\geq\mathbf{1}_{\{\rho=0\}}\beta^n_{+}$. As the very last step, we see that $\mathbf{1}_{\rho=0}\beta^n_{+}=\mathbf{1}_{\rho=0}\beta^n$. Indeed, this follows from an application of \cite[Lemma 2.10]{DDNL18a} for any $x_0\in{\{\rho=0}\},u:=\rho=g$ and $q:=g$, where $g$ is a local potential of $\beta$ near $x_0$.
\end{proof}
As an important role in studying Monge-Amp\`ere equations with prescibed singularity profile in the sense of Darvas-Di Nezza-Lu, the so-clalled rooftop envelope is defined by
$$P[\psi](\phi):=(\lim_{C\to +\infty}P(\psi+C,\phi))^{*}$$ where $\psi,\phi\in\mbox{PSH}(X,\beta)$.
If $\phi=\rho$, we denote $P[\psi](\rho)$ by $P[\psi]$ for simiplicity. This envelope has good properties which can be reformulated from Darvas-Di Nezza-Lu's papers. By the proof of \cite[Lemma 2.9, Theorem 3.6]{DDNL23} and Theorem \ref{regularity envelop}, we have
\begin{cor}\label{more general regularity of envelope}
Suppose $\phi,\psi,P(\phi,\psi)\in\mbox{PSH}(X,\beta)$. Then
$$\mbox{MA}(P(\phi,\psi))\leq\mathbf{1}_{\{P(\phi,\psi)=\phi\}}\mbox{MA}(\phi)+\mathbf{1}_{\{P(\phi,\psi)=\psi\}}\mbox{MA}(\psi).$$
In particular,
$$\mbox{MA}(P[\psi](\phi))\leq\mathbf{1}_{\{P[\psi](\phi)=\phi\}}\mbox{MA}(\phi)~and~\mbox{MA}(P[\psi])\leq\mathbf{1}_{\{P[\psi]=0\}}\beta^n.$$
\end{cor}
Using the result in \cite{LWZ23} and Corollary \ref{more general regularity of envelope}, we can prove a non-collapsing property for the class of potentials with the same singularity type as $\phi$, when $\int_X \mbox{MA}(\phi)>0$. See the proof of \cite[Corollary 3.9]{DDNL18b}.
\begin{cor}\label{ab}
Assume that $\phi\in\mbox{PSH}(X,\beta)$ is such that $\int_X \mbox{MA}(\phi)>0$. If $U$ is a Borel subset of $X$ with positive Lebesgue measure, then there exists $\psi\in\mbox{PSH}(X,\beta)$ having the same singularity type as $\phi$ such that $\mbox{MA}(\psi)(U)>0$.
\end{cor}
Using the proof of \cite[Theorem 3.14]{DDNL23}, Corollary \ref{more general regularity of envelope} has an important application:
\begin{cor}
Assume that $\phi\in\mbox{PSH}(X,\beta)$ and $\int_X \mbox{MA}(\phi)>0$. Then 
$$P[\phi]=\sup_{v\in F_{\phi}}v,$$
where $F_{\phi}:=\{v\in\mbox{PSH}(X,\beta):\phi\leq v\leq0~and~\int_X \mbox{MA}(v)=\int_X \mbox{MA}(\phi)\}$. In particular, $P[\phi]=P[P[\phi]]$.
\end{cor}

\section{Solution of Monge-Amp\`ere equations with prescribed singularities}

Now we turn to define the relative pluripotential concepts in line with Darvas-Di Nezza-Lu. 
\begin{defn}\label{def of relative full mass class}
Given a quasi-psh function $\phi\in\mbox{PSH}(X,\beta)$, the relative full mass class $\mathcal{E}(X,\beta,\phi)$ is the set of all $\beta$-psh functions $u$ satisfying that $u$ is more singular than $\phi$ and $\int_X(\beta+dd^c u)^n=\int_X (\beta+dd^c \phi)^n$. 
\end{defn}
By Lemma \ref{same sing same mass} we know that $\mathcal{E}(X,\beta,\phi)$ contains all $\beta$-psh functions $u$ which have the same singularity type with $\phi$. By the proof of \cite[Remark 3.8]{DDNL23}, $\mathcal{E}(X,\beta,\phi)$ is stable under taking the maximum.
\begin{defn}\label{def of model potential}
A model potential is a $\beta$-psh function $\phi$ such that $P[\phi]=\phi$ and $\int_X \mbox{MA}(\phi)>0$.
\end{defn}

It is easy to see that all bounded $\beta$-psh functions are model potentials. After we prove the existence of Monge-Amp\`ere equation with prescribed singularity type, we can prove that all $\beta$-psh functions with analytic singularities are also model potentials, see \cite[Proposition 4.36]{DDNL18b}\cite[Proposition 5.23]{DDNL23}. For simplicity we always assume that $\phi$ is a model potential from now on.

As Darvas-Di Nezza-Lu did, we introduce the relative Monge-Ampère capacity of a Borel set $B\subset X$:
$$\mbox{Cap}_{\phi}(B):=\sup\{\int_B (\beta+dd^c \psi)^n,\psi\in\mbox{PSH}(X,\beta),\phi\leq\psi\leq\phi+1\}.$$
Using the proof of \cite[Lemma 4.3]{DDNL23} and Theorem \ref{local polar=global polar}, we have that
\begin{thm}\label{21}
Let $B\subset X$ be a Borel subset. Then $\mbox{Cap}_{\phi}(B)=0$ iff $B$ is pluripolar.
\end{thm}
\begin{proof}
By Theorem \ref{local polar=global polar} a Borel set $E\subset X$ is pluripolar iff $\mbox{Cap}_{\beta}(E)=0$. If $B$ is pluripolar then $\mbox{Cap}_{\phi}(B)=0$ by definition. Conversly, assume that $\mbox{Cap}_{\phi}(B)=0$.If $B$ is non-pluripolar then $\mbox{Cap}_{\beta}(B)>0$. Since $\mbox{Cap}_{\beta}$ is inner regular (\cite[Lemma 4.2]{DDNL23}), there exists a compact subset $K$ of $B$ such that $\mbox{Cap}_{\beta}(K)>0$. In particular $K$ is non-pluripolar, hence the global extremal function $V_{K,\rho}$ is bounded from above. Then see the proof of \cite[Lemma 4.3]{DDNL18b}.
\end{proof}
Now we reformulate some basic properties in \cite{DDNL18b}\cite{DDNL21} under our setting. We denote the set of all $\beta$-psh functions which are more singular than $\phi$ by $PSH(X,\beta,\phi)$. $\phi$ is said to have small unbounded locus if it is locally bounded outside a closed complete pluripolar subset.
Using Theorem \ref{convergence thm}, Corollary \ref{more general regularity of envelope}  and the proof of \cite[\S 3.2]{DDNL21} we can prove that:

\begin{prop}
Let $E$ be a Borel subset of $X$ and $h_{E,\phi}$ be the relative extremal function defined below:
$$h_{E,\phi}:=\sup\{u\in\mbox{PSH}(X,\beta,\phi)|u\leq\phi-1~on~E;u\leq0~on~X\}.$$ 
Then $h_{E,\phi}^{*}$ is a $\beta$-psh function such that $\phi-1\leq h_{E,\phi}^{*}\leq\phi$ and $\mbox{MA}(h_{E,\phi}^{*})$ vanished on $\{h_{E,\phi}^{*}<0\}\backslash\bar{E}$.
\end{prop}
\begin{prop}
If $E$ is a non-pluripolar subset of $X$ then $\mbox{MA}(V_{E,\phi}^{*})$ vanishes in $X\backslash\bar{E}$, where $V_{E,\phi}^{*}$ is defined by
$$V_{E,\phi}:=\sup\{u\in\mbox{PSH}(X,\beta,\phi)|u\leq\phi~on~E\}.$$
\end{prop}
\begin{lem}
Assume that $K$ is a compact subset of $X$ and $\mbox{Cap}_{\phi}(K)>0$. Then we have
$$1\leq(\frac{\int_X \mbox{MA}(\phi)}{\mbox{Cap}_{\phi}(K)})^{1/n}\leq\max(1,M_{\phi}(K)),$$
where $M_{\phi}(K):=\sup_X V_{E,\phi}^{*}$.
\end{lem}
Using the lemma above and following the line in \cite[Proposition 4.30]{DDNL21} we have
\begin{prop}
Let $f\in L^p(X,\omega^n),p>1$ with $f>0$. Then there exists a constant $C>0$ depending only on $\beta,\omega,X,n,||f||_{L^p}$ such that
$$\int_E f\omega^n\leq\frac{C}{\int_X \mbox{MA}(\phi)}\cdot\mbox{Cap}_{\phi}(E)^2$$
for every Borel subsets $E\subset X$.
\end{prop}

\subsection{Variational method}
In this part we assume that $\phi$ has small unbounded locus to ensure that the integration-by-parts formula is valid. We follow the line of \cite[\S 4]{DDNL18b}. To develope the variational approach, we define the relative Monge-Amp\`ere energy of $u\in\mathcal{E}(X,\beta,\phi)$ such that $u$ has the same singularity type as $\phi$ by
$$E_{\phi}(u):=\frac{1}{n+1}\sum_{k=0}^{n}\int_X (u-\phi)(\beta+dd^c u)^k\wedge(\beta+dd^c \phi)^{n-k}.$$
For arbitrary $u\in\mbox{PSH}(X,\beta,\phi)$ we define its Monge-Amp\`ere energy by
$$E_{\phi}(u):=\inf\{E_{\phi}(v)|v\in\mathcal{E}(X,\beta,\phi),v~has~the~same~singularity~type~as~\phi~and~u\leq v\}.$$
Note that these definitions coincide for $u\in\mathcal{E}(X,\beta,\phi)$ having the same singularity type as $\phi$, see \cite[Lemma 4.11]{DDNL18b}.
The energy operator enjoys basic properties which goes back to \cite{BEGZ10}, see \cite[\S 4.2]{DDNL18b}.
\begin{defn}\label{defn of finite energy class}
$\mathcal{E}^1(X,\beta,\phi)$ is defined to be the set of all $u\in\mbox{PSH}(X,\beta,\phi)$ such that $E_{\phi}(u)$ is finite. 
\end{defn}
By \cite[Lemma 4.12, Theorem 4.10]{DDNL21} we know that $E_{\phi}$ is non-decreasing and $\mathcal{E}^1(X,\beta,\phi)$ is stable under taking the maximum. For simiplicity, $\mathcal{E}(X,\beta):=\mathcal{E}(X,\beta,0)$ and $\mathcal{E}^1(X,\beta):=\mathcal{E}^1(X,\beta,0)$ if $\phi=0$.
Now we can describe the first order variation of $E_{\phi}$:
\begin{prop}\label{20}
Let $u\in\mathcal{E}^1(X,\beta,\phi)$ and $\chi$ be a continuous function on $X$. For each $t>0$ we set $u_t:=P(u+t\chi)$. Then $u_t\in\mathcal{E}^1(X,\beta,\phi),t\mapsto E_{\phi}(u_t)$ is differentiable and its derivative is given by
$$\frac{d}{dt}E_{\phi}(u_t)=\int_X \chi\mbox{MA}(u_t),t\in\mathbf{R}.$$
\end{prop}
\begin{proof}
Note that Proposition \ref{mass concentr of envelope} implies that $\mbox{MA}(u_t)$ is supported on $\{u_t=u+t\chi\}$. Follow the line of \cite[Lemma 4.20]{DDNL18b}.
\end{proof}
We introduce the following functionals on $\mathcal{E}^1(X,\beta,\phi)$:
$$F_{\lambda}(u):=F_{\lambda,\mu}(u):=\mathrm{E}_{\phi}(u)-L_{\lambda,\mu}(u),u\in\mathcal{E}^{1}(X,\beta,\phi),$$ where $L_{\lambda,\mu}(u):=\frac{1}{\lambda}\int_{X}e^{\lambda u}d\mu $ if $\lambda>0$ and $L_{\mu}(u):=L_{0,\mu}(u):=\int_{X}(u-\phi)d\mu$. For each constant $A\geq1$ we set ${\mathcal{M}}_{A}$ denote the set of all probability measures $\mu$ on X such that
$$\mu(E)\leq A\cdot\text{Cap}_\phi(E),\text{for all Borel subsets}~E\subset X.$$
By the proof of \cite[Lemma 4.17]{DDNL18b}, $\mathcal{M}_A$ is a compact convex subset of the set of probability measures on $X$.
\begin{thm}
Assume that $L_{\lambda,\mu}$ is finite on $\mathcal{E}^1(X,\beta,\phi)$ and $u\in\mathcal{E}^1(X,\beta,\phi)$ maximizes $F_{\lambda,\mu}$ on $\mathcal{E}^1(X,\beta,\phi)$. Then $u$ solves the Monge-Amp\`ere equation
$$(\beta+dd^c u)^n=e^{\lambda u}\mu,\lambda\geq0.$$

\end{thm}
\begin{proof}
We may assume that $\lambda>0$. Let $\chi$ be an arbitrary continuous function on $X$ and set $u_t:=P(u+t\chi)$. By Proposition \ref{20} $u_t\in\mathcal{E}^1(X,\beta,\phi)$ for all $t$ and $$g(t):={E}_\phi(u_t)-L_{\lambda,\mu}(u+t\chi)$$
is differentiable on $\mathbf{R}$ whose derivative is given by $g^{\prime}(t)=\int_{X}\chi\mbox{MA}(u_t)-\int_{X}\chi e^{\lambda(u+t\chi)}d\mu.$ Since $u_{t}\leq u+t\chi$, $g(t)\leq F_{\lambda,\mu}(u_{t})\leq\operatorname*{sup}_{\mathcal{E}^{1}(X,\beta,\phi)}{F}_{\lambda,\mu}=F(u)=g(0),$ which implies $g$ attains a maximum at 0. Hence $g^{\prime}(0)=0$. Since $\chi$ is arbitrary it follows that $\mbox{MA}(u)=e^{\lambda u}\mu$. If $\lambda=0$, the result follows from a similar argument.
\end{proof}
As an important application of the theorem above (see \cite[Theorem 4.25]{DDNL18b}), we have
\begin{thm}
Assume that $\mu\in\mathcal{M}_A$ for some $A\geq1$. Then there exists $u\in\mathcal{E}^1(X,\beta,\phi)$ such that $(\beta+dd^c u)^n=\mu$.
\end{thm}

Now we can prove our first main theorem. 
\begin{thm}
Assume that $\nu$ is a positive non-pluripolar measure on $X$ such that $\mu(X)=\int_X \mbox{MA}(\phi)$. Then there exists $u\in\mathcal{E}(X,\beta,\phi)$ (unique up to a constant) such that $(\beta+dd^c u)^n=\mu$.
\end{thm}
\begin{proof}
By Theorem \ref{21}  we know that for every Borel subset $E$ is pluripolar if and only if $\mbox{Cap}_{\phi}(E)=0$. By the proof of \cite[Lemma 4.26]{DDNL18b} (see also \cite{Ce98}), $\mu=f\nu$ where $\nu\in\mathcal{M}_1$ and $0\leq f\in L^1(X,\nu)$. The rest of the argument is identical with \cite[Theorem 4.28]{DDNL18b}.
\end{proof}
\begin{rem}\label{lambda>0}
By a similar argument to \cite[Theorem 4.23]{DDNL18b} there exists a unique $\varphi\in\mathcal{E}^1(X,\beta,\phi)$ such that $(\beta+dd^c \varphi)^n=e^{\lambda \varphi}\mu$ with $\lambda>0$. On the other hand, we cannot use the proof of \cite[Theorem 5.17]{DDNL23} to prove our main theorem above. This is because although it is easy to see that the  Monge-Amp\`ere  capacity $\mbox{Cap}_{\beta}$ can be dominated by the Bedford-Taylor capacity $\mbox{Cap}_{BT}$, it is not known for the converse. 
\end{rem}
\begin{rem}
By a similar argument of \cite[Corollary 4.33]{DDNL18b}, the regularity theorem is also valid under our setting.
\end{rem}
\begin{rem}
The theorem above generalizes the result of \cite{LWZ23} in the relative context. On the other hand, if we focus on the special case that $\phi=\rho$, we can use the supersolution techniques under our setting, which means that we can relax the condition of model potential $\phi$ that it has small unbounded locus. One only need to follow the line of \cite[\S4]{DDNL21}.
\end{rem}
\subsection{Weighted energy class}
Now we want to prove our second main theorem. We adapt the proof from \cite{DDNL18a}. In what follows, $\chi$ is denoted as a continuous increasing function from $[0,+\infty)$ to $[0,+\infty)$ such that $\chi(0)=0,\chi(\infty)=\infty$ and satisfies the following condition:
$$\forall t>0,\forall\lambda>1,\chi(\lambda t)\leq\lambda^{M}\chi(t)$$ where $M\geq1$ is a fixed constant. It is easy to see that for every $t,s>0$, $\chi(t+s)\leq 2^{M}\max(\chi(t),\chi(s))$.

Fix $\phi$ a model potential and let $\mathcal{E}_{\chi}(X,\beta,\phi)$ denote the set of all $u\in\mathcal{E}(X,\beta,\phi)$ such that 
$$E_{\chi}(u,\phi):=\int_X \chi(|u-\phi|)\mbox{MA}(u)<\infty.$$
Now we need to establish some important properties of the weighted energy class defined above.
\begin{thm}\label{25}
If $\varphi,\psi\in\mathcal{E}(X,\beta,\phi)$, then $P(\varphi,\psi)\in\mathcal{E}(X,\beta,\phi)$. In particular, if $\varphi,\psi\in\mathcal{E}(X,\beta,\phi)$, then $P(\varphi,\psi)\in\mathcal{E}(X,\beta)$.
\end{thm}
\begin{proof}
We may assume that $\varphi,\psi\leq0$. Let $\varphi_{j}:=\operatorname*{max}(\varphi,\phi-j)$ and $\psi_j:=\max(\psi,\phi-j)$. For each $j>0$, it follows from Lemma \ref{24} below that there exists a unique $u_j\in\mbox{PSH}(X,\beta,\phi)$ having the same singularity type with $\phi$ such that $$
(\beta+dd^c u_j)^n=e^{u_j-\varphi_j}(\beta+dd^c \varphi_j)^n+e^{u_j-\psi_j}(\beta+dd^c\psi_j)^n.$$
It follows from the proof of \cite[Lemma 2.5]{DDNL18a} (see also domination principle \cite[Theorem 3.12]{DDNL23} and comparison principle \cite[Corollary 3.23]{DDNL23}) that $u_j\leq\min(\varphi_j,\psi_j)$. Hence $u_j\leq P(\varphi_j,\psi_j)$. Next we claim that $$\sup_j\int_X \chi(\phi-u_j)(\beta+dd^c u_j)^n<+\infty.$$ To prove the claim, in view of the equation above, it is enough to prove that 
$$\sup_j\int_X \chi(\phi-u_j)e^{u_j-\varphi_j}(\beta+dd^c \varphi_j)^n<+\infty.$$ By definition of $\chi$ it suffices to check that 
$$\sup_j\int_X \chi(\phi-\varphi_j)e^{u_j-\varphi_j}(\beta+dd^c \varphi_j)^n<+\infty.$$ By the proof of \cite[Lemma 6.5]{DDNL23} and $u_j\leq\varphi_j$ and $\varphi\in\mathcal{E}_{\chi}(X,\beta,\phi)$, the claim holds.

Since $\chi(+\infty)=+\infty$, the claim implies that $\sup_X u_j$ is uniformly bounded. It thus follows from the proof of \cite[Lemma 6.7]{DDNL23} some subsequence of $u_j$ converges in $L^1(X,\omega)$ to some $u\in\mathcal{E}_{\chi}(X,\beta,\phi)$. Since $u_j\leq P(\varphi_j,\psi_j)$, we have $u\leq P(\varphi,\psi)$. By the proof of \cite[Lemma 2.5]{DDNL23} and \cite[Corollary 3.19]{DDNL23} the first statement holds. The second statement follows from Lemma \ref{23} and the first statement.
\end{proof}
The following lemma is essentially known in \cite{BEGZ10}:
\begin{lem}\label{23}
Let $\varphi$ be a $\beta$-psh function. Then it has full relative Monge–Ampère mass if and only if $E_{\chi}(\varphi,\phi)<\infty$ for some weight function $\chi$.
\end{lem}
\begin{proof}
Since $\phi$ is a model potential, we know $\varphi\leq\phi$. Let $\varphi^{(k)}:=\max\{\varphi,\phi-k\}$. It suffices to show that 
$$ m_k:=\int_{\{\varphi\leqslant\phi-k\}}\mathrm{MA}(\varphi^{(k)})$$
tends to 0 iff
$$ \sup_k\int_X\chi(\phi-\varphi^{(k)})\mathrm{MA}(\varphi^{(k)})<\infty $$
for some weight $\chi$. We have 
$$\begin{aligned}
\int_{X}\chi(\phi-\varphi^{(k)})\mathrm{MA}(\varphi^{(k)}) &=\chi(k)m_{k}+\int_{\{\varphi>\phi-k\}}\chi(\phi\varphi)\mathrm{MA}(\varphi) \\
&\leq\chi(k)m_{k}+\int_{0}^{+\infty}\mbox{MA}(\varphi)(\{\varphi\leq\rho-\chi^{-1}(t)\})dt\\
&\leq\chi(k)m_{k}+\int_{0}^{+\infty}\mbox{MA}(\varphi)(\{\varphi\leq\rho-t\})\chi^{\prime}(t)dt.
\end{aligned}$$
Since $\mbox{MA}(\varphi)$ puts no mass on the pluripolar set, we can construct a piecewise linear function $\chi$ such that $\chi(k)\leq\frac{1}{m_k}$ and $\chi^{\prime}(k)\leq\frac{1}{k^2}$ when $k$ is near $+\infty$ and $\chi^{\prime}(k)=1$ when $k$ is near 0. Then the result holds.

\end{proof}
\begin{lem}\label{24}
Assume that $u,v\in\mbox{PSH}(X,\beta)$ having the same singularity type with $\phi$. Then there exists a unique $\varphi\in\mbox{PSH}(X,\beta)$ having the same singularity type with $\phi$ and satisfying
$$(\beta+dd^c \varphi)^n=e^{\varphi-u}(\beta+dd^c u)^n+e^{\varphi-v}(\beta+dd^c v)^n.$$
\end{lem}
\begin{proof}
To prove the unqiueness we need to use the proof of \cite[Lemma 5.15]{DDNL23}. To prove the existence we need to  approximate $u$ by $u_j:=\max(u,\rho-j)$ and use the Remark \ref{lambda>0}. Then use the proof of \cite[Lemma 2.14]{DDNL18a}.
\end{proof}
By a similar argument to \cite[Theorem 2.17]{DDNL18a}, it follows from Proposition \ref{mass concentr of envelope} and Theorem \ref{25} that
\begin{thm}
Assume that $\psi,\varphi\in\mathcal{E}(X,\beta,\phi)$. Then $P[\varphi](\psi):=(\lim_{C\to +\infty}P(\varphi+C,\psi))^{*}=\psi$.
\end{thm}
Now we can prove our second main theorem.
\begin{thm}\label{lelong0}
Assume that $\varphi\in\mathcal{E}(X,\beta,\phi)$. Then $\forall x\in X$, 
$$\nu(\varphi,x)=\nu(\phi,x).$$
\end{thm}
\begin{proof}
Since $\phi$ is a model potential, $\varphi\leq\phi$ and hence $\nu(\varphi,x)\geq\nu(\phi,x)$ for every $x\in X$. We will argue by contradiction. Suppose that $\nu(\varphi,x)>\nu(\phi,x)$ for some $x$.  By the definition of Lelong number, there exists a coordinate ball $\mathbb{B}$ centered at $x$, such that
\begin{align*}
\varphi(z)\leq \gamma \log |z|+C\ \mbox{in}\ \mathbb{B}.
\end{align*}
It follows from Theorem \ref{25} that $P[\varphi](\phi)=\phi$. Let $g$ be the local potential of $\beta$ such that $g+\phi,g+\varphi\leq0$ in $\mathbb{B}$. In view of the definition of $P[\varphi](\phi)$, it holds that
\begin{align*}
g+\phi=g+P_{[\beta,\varphi]}(\phi)\leq\sup\{v\in\mbox{PSH}(\mathbb{B}):v\leq0, v\leq \gamma\log|z|+O(1)\}.
\end{align*}
The RHS is the pluricomplex Green function $G_{\mathbb{B}}(z,0)$ with logarithmic pole at $0$ of order $\gamma$, which is equal to $\gamma\log|z|$ by definition. Comparing the Lelong number of both sides, it makes a contradiction.
\end{proof}
\begin{rem}
If we assume that $X$ is K\"ahler and $\phi=0$, $\beta$ is semi-positive, the theorem was proved by Darvas-Di Nezza-Lu. In fact, their result also holds on compact complex manifold $Y$ if there exists a modification $f:X\to Y$, i.e. of Fujiki class. Assume that $\alpha$ is a semi-positive and big closed form on $Y$, we pick a function $u\in PSH(Y,\alpha)$. Taking a local potential $v$ of $\alpha$ at any point $p$ in the exceptional locus in $Y$, we have $u+v$ is psh and $f^{*}(u)\in PSH(X,f^{*}\alpha)$ where $f^{*}\alpha$ is a semi-positive and big form on $X$. By Zariski's main theorem we know that the fibre $f^{-1}(p)$ is connected and by compactness $f^{*}u$ is constant in $f^{-1}(p)$. Theorefore $\nu(f^{*}u)=\nu(u)=0$ by Darvas-Di Nezzaa-Lu's theorem.
\end{rem}

\begin{rem}
Following the same trick as above, same conclusions hold when $\omega$ is K\"ahler and $\beta=\theta$ a smooth representative of a big class on $X$. This gives a relative version of generalization on the results in \cite{DDNL18a}.
\end{rem}
At the end of this section, we give a simple application of Theorem \ref{lelong0} under the assumption that $\phi=\rho$. See also \cite{BBGZ13}. 
We need the following lemma first:
\begin{lem}\label{lem:continuity of L minus}
The map $\mathcal{E}^1(X,\beta)\rightarrow L^1(X,\omega^n),\varphi\mapsto e^{-\varphi}$ is continuous.
\end{lem}
\begin{proof}
Thanks to Theorem \ref{lelong0} and Skoda's integrable theorem, we know that for $\varphi\in\mathcal{E}(X,\beta)$, $e^{-\varphi}$ is  $L^1$ integrable with respect to Lebesgue measure. Then, by \cite[Theorem 1.1]{GuLiZh22}, we know that if $\varphi_j\rightarrow\varphi$ in $\mathcal{E}^1(X,\beta)$, then $e^{-\varphi_j}\rightarrow e^{-\varphi}$ in $L^1(X,\omega^n)$, which completes the proof.
\end{proof}
Set
\begin{align*}
L_{-}(\varphi):=-\log \int_Xe^{-\varphi}\omega^n
\end{align*}
and $F_{-}(\varphi)=-E_{\rho}(\varphi)-L_{-}(\varphi)$, which is well-known as Ding functional.
Our main theorem about the application of Theorem \ref{lelong0} is the following:
\begin{thm}\label{thm:MA minus}
Assume there exists $\varphi\in\mathcal{E}^1(X,\beta)$ such that $F_{-}(\varphi)=\sup_{\mathcal{E}^1(X,\beta)}F_{-}$, then $\varphi$ is bounded on $X$ and for some $c\in\mathbb{R}$,
\begin{align*}
(\beta+dd^c\varphi)^n=e^{-(\varphi+c)}\omega^n.
\end{align*}
\end{thm}
\begin{rem}
In view of Lemma \ref{lem:continuity of L minus}, we know that $F_{-}$ is upper-semi-continuous on $\mathcal{E}^1(X,\beta)$. Hence the properness (in the sense of \cite{BBGZ13}) of $F_{-}$ will lead to the existence of solution to the Monge-Amp\`ere equation.
\end{rem}
\noindent{\bf Proof of Theorem \ref{thm:MA minus}:}
Take $\chi\in C(X)$, since $\varphi$ achieves the supremum of $F_{-}$, we get that
\begin{align*}
E_{\rho}(P(\varphi+t\chi))+\log \int_Xe^{-(\varphi+t\chi)}&\leq E_{\rho}(P(\varphi+t\chi))+\log \int_X e^{-P(\varphi+t\chi)}\\
&\leq E_{\rho}(\varphi)+\log \int_Xe^{-\varphi}.
\end{align*}
Thus the function $g:t\mapsto E_{\rho}(P(\varphi+t\chi))-L_{-}(\varphi+t\chi)$ achieves its supremum at 0. By Proposition \ref{20}, we know that $g$ is differentiable at $t=0$, hence $g'(0)=0$. This leads to 
\begin{align*}
0=g'(0)=\int_X\chi(\beta+dd^c\varphi)^n-
\frac{\int_X\chi e^{-\varphi}\omega^n}{\int_Xe^{-\varphi}\omega^n}.
\end{align*}
Since $\chi$ is arbitrary, we get that for $c=\log \int_Xe^{-\varphi}\omega^n$,
\begin{align*}
(\beta+dd^c\varphi)^n=e^{-(\varphi+c)}\omega^n.
\end{align*}
Note that $\varphi$ has zero Lelong number on $X$, this gives that $e^{-\varphi}$ is actually $L^p$-integrable for any $p>1$. It then follows from $L^{\infty}$-estimate in \cite[\S 3.4]{LWZ23} that $\varphi$ is bounded on $X$. This completes the proof.

\end{document}